\documentclass[10pt,a4paper,twoside]{amsart}

\usepackage{amsgen}
\usepackage{amsmath}
\usepackage{amstext}
\usepackage{amsbsy}
\usepackage{amsopn}
\usepackage{amsfonts}
\usepackage{amsthm}
\usepackage{amscd}
\usepackage{amsxtra}
\usepackage{eucal}
\usepackage{amscd}
\usepackage{amsfonts}
\usepackage{epsf}
\usepackage[all]{xy}
\usepackage{amssymb}
\usepackage{pb-diagram}
\DeclareMathOperator{\Der}{Der}
\DeclareMathOperator{\Image}{Im}
\DeclareMathOperator{\Ker}{Ker}

\newtheorem{theorem}{Theorem}[section]
\newtheorem*{theorem*}{Theorem}
\newtheorem{lemma}[theorem]{Lemma}
\newtheorem{proposition}[theorem]{Proposition}
\newtheorem{corollary}[theorem]{Corollary}
\theoremstyle{definition}
\newtheorem{definition}[theorem]{Definition}

\theoremstyle{remark}
\newtheorem{remark}[theorem]{Remark}
\newtheorem*{nota}{Notation}
\numberwithin{equation}{section}

\newcommand{\ord}{\operatorname{ord}}
\newcommand{\initform}{\operatorname{in}}
\newcommand{\Spec}{\operatorname{Spec}}
\newcommand{\Proj}{\operatorname{Proj}}
\newcommand{\Hom}{\operatorname{Hom}}

\newcounter{temp}

\makeatletter
\def\square{\RIfM@\bgroup\else$\bgroup\aftergroup$\fi
\vcenter{\hrule\hbox{\vrule\@height.6em\kern.6em\vrule}\hrule}\egroup}
\makeatother

\begin{document}
\bibliographystyle{amsalpha}

\title[Infinitesimal deformations and obstructions]{On infinitesimal
deformations and obstructions for rational surface singularities}

\author[Christophersen]{Jan Arthur Christophersen}
\address{Department of Mathematics \\
P.O. 1053 Blindern\\
N-0316 OSLO, Norway}
\email{christop@math.uio.no}

\author[Gustavsen]{Trond St\o len Gustavsen}
\address{Department of Mathematics \\
P.O. 1053 Blindern\\
N-0316 OSLO, Norway}
\email{stolen@math.uio.no}

\date{\today}

\maketitle
\section*{Introduction}

The purpose of this paper is to prove dimension formulas for $T^1$
and $T^2$ for rational surface singularities. These modules play an
important role in the deformation theory of isolated singularities in
analytic and algebraic geometry. The first may be identified as the
Zariski tangent space of the versal deformation of the singularity;
i.e. it is the space of infinitesimal deformations. The second
contains the obstruction space -- in all known cases it is the
whole obstruction space for rational surface singularities.

The dimension formulas for $T^1_X$ and $T^2_X$ relate these
dimensions to similar invariants on the blow up, $\widehat{X}$, of
$X$.
An important result of Tjurina, which we state below
(Theorem~\ref{Tyu}), shows that the minimal resolution
$\widetilde{X}$
may be gotten by a series of blow-ups. Thus, in principle, the
formulas allow one to compute these dimensions via blowing up. In
fact, the nature of the formulas allows one in many cases to compute
these dimensions from the graph.

Computing $T^1$ and $T^2$ for rational surface singularities has a
history which we briefly recall. (The terms involved here are
explained in Section~\ref{prelim}.) Of course for the rational double
points $T^2=0$ and $T^1$ is easily computed. We will from now on
assume that singularities are not hypersurfaces; i.e. the embedding
dimension $e$ is not $3$. In the 80's much work was done in Hamburg
on computing $T^1$ for quotient surface singularities, a sub-set of
the rationals, and the general form turned out to be $\dim
T^1_X=(e-4) + \dim H^1(\widetilde{X}, \Theta_{\widetilde{X}})$
(\cite{bkr:inf}).
Behnke and Kn\"{o}rrer (\cite{bk:inf}) where able to prove the same
formula for a larger, but still very restricted class of rational
surface singularities. In 1987, J. Arndt and the first author proved
independently that for a cyclic quotient singularity $\dim
T^2_X=(e-2)(e-4)$. Later, using hypersurface sections, Behnke and the
first author proved this formula for rational surface singularities
with reduced fundamental cycle and $T^2_{\widehat{X}}=0$ and for all
quotient singularities (\cite{bc:hypr}). Finally, de Jong and van
Straten
(\cite{js:rss}), gave the correct formulas for {\it all} rational
surface singularities with reduced fundamental cycle.

Let $I$ be an index set for all singularities $X_{\nu}$ (including
$X$ itself) that appear in the process of resolving a rational
surface singularity with reduced fundamental cycle by blowing up
points. Let $e(\nu)$ be the embedding dimension of $X_{\nu}$ and set
$I_4=\{{\nu}\in I:e(\nu)\ge 4\}$; i.e. the indices of
non-hypersurface
singularities. What de Jong and van Straten proved (\cite[Theorem
3.16]{js:rss}) was that
\begin{align*}
\dim T^1_X&=\sum_{{\nu}\in I_4}(e(\nu)-4) +
        \dim H^1(\widetilde{X}, \Theta_{\widetilde{X}})\\
\dim T^2_X&=\sum_{{\nu}\in I_4}(e(\nu)-2)(e(\nu)-4)\, .
\end{align*}
On $\widehat{X}$ we have the groups $T^1_{\widehat{X}}$ and
$T^2_{\widehat{X}}$
(see Section~\ref{cotangent}). Since rational surface singularities
are absolutely isolated, $\dim
T^1_{\widehat{X}}=\sum_{p\in\widehat{X}}\dim
T^1_{\widehat{X},p}+ \dim H^1(\widehat{X}, \Theta_{\widehat{X}})$ and
$\dim T^2_{\widehat{X}}=\sum_{p\in\widehat{X}}\dim
T^2_{\widehat{X},p}$.
It is well known that for a rational double point $\dim T^1_X=
\dim H^1(\widetilde{X}, \Theta_{\widetilde{X}})$. Thus, using the
Leray
spectral sequence for $p:\widetilde{X}\to\widehat{X}$ and that
$p_{\star}\Theta_{\widetilde{X}}\simeq \Theta_{\widehat{X}}$ we see
that the
de Jong--van Straten result is equivalent to saying that for all
rational surface singularities with reduced fundamental cycle
\begin{align*}
\dim T^1_X&=(e-4) +\dim T^1_{\widehat{X}}\\
\dim T^2_X&=(e-2)(e-4)+\dim T^2_{\widehat{X}}\, .
\end{align*}

The results in this paper originated from a wish to find a direct
relationship between the $T^i$ and blowing up for rational
singularities. This is described in Section~\ref{cotangent}. This
allows us to compute the $T^i$ in terms of the cohomology of certain
sheaves on $\widehat{X}$. What we get (Theorem~\ref{t1} and
Theorem~\ref{t2}) is that for all rational surface singularities
(with $e\ge 4$)
\begin{align*}
\dim T^1_X&=(e-4) +\dim T^1_{\widehat{X}}+c(X)\\
\dim T^2_X&=(e-2)(e-4)+\dim T^2_{\widehat{X}}+c(X)
\end{align*}
where $c(X)$ is the dimension of the $H^1$ of a certain sheaf (in
fact several) on $\widehat{X}$ (Definition~\ref{cor}). We give some
partial results on $c(X)$ in Section~\ref{Cor_term}, in particular we
show that $c(X)=0$ when the fundamental cycle is reduced, reproving
the de
Jong--van Straten result.

\subsection*{Acknowledgment} The results in this paper have developed
through continuous discussions with Kurt Behnke, Theo de Jong, Jan
Stevens and Duco van Straten.  Jan Stevens pointed out a serious
mistake in the first version of this paper.  We are grateful to Olav
Arnfinn Laudal for patiently answering questions and explaining the
cohomology theory involved and to Kristian Ranestad for helping us out
with the geometry.

\section{Preliminaries}\label{prelim}

\subsection{Results on rational singularities}
The singularities we study are algebraic over $\mathbb C$, i.e. of
the form $X=\Spec A$ where $A=P/I$ and $P$ is a regular local
$\mathbb C$ algebra essentially of finite type.
A normal surface singularity $X$ with minimal resolution
$f:\widetilde{X}\to X$ is rational if $H^1(\widetilde{X},{\mathcal
O}_{\widetilde{X}})=0$ (\cite{ar:rat}). The exceptional divisor
$E\subset
\widetilde{X}$ is a union of irreducible components $E_i\simeq
{\mathbb
P}^1$. There is a {\it fundamental cycle} $Z$, supported on $E$,
defined by ${\mathfrak m}{\mathcal O}_{\widetilde{X}}$. Here
${\mathfrak
m}$ is the maximal ideal in ${\mathcal O}_X$. This divisor may be
constructed as the unique smallest positive divisor $Z=\sum r_iE_i$
satisfying $Z\cdot E_i\le 0$ for all irreducible components $E_i$.
The embedding dimension of $X$, $e=\dim_{\mathbb C}{\mathfrak
m}/{\mathfrak m}^2$, equals $-Z^2+1$ and the multiplicity
$m(X)=e-1=-Z^2$.

There are three theorems on rational surface singularities which are
essential for our results. We collect them and partially rephrase
them here. The first is a result from \cite{tj:abs} which shows how
the
blow up $\widehat{X}$ may be obtained from $\widetilde{X}$.
\begin{theorem}[Tjurina]\label{Tyu} If $X$ is a rational surface
singularity, then the blow up of $X$ is isomorphic to the surface
obtained from $\widetilde{X}$ by contracting all components $E_i$
with
$Z\cdot E_i=0$.
\end{theorem}

In \cite{wa:equ} we find the basic algebraic property of rational
surface singularities that we will need.
\setcounter{temp}{0}
\begin{theorem}[Wahl] \label{Wah}
Let $A=P/I$ be the local ring of a rational surface singularity,
where $P$ is a regular local algebra of dimension $e$ over an
algebraically closed field $k$. Let $\bar{P}$ and $\bar{A}$ be the
associated graded rings with respect to the maximal ideal. Then there
exist minimal projective resolutions:
\begin{gather*}
0\xrightarrow{} P^{b_{e-2}}\xrightarrow{\phi_{e-2}}\dots
\xrightarrow{} P^{b_2}\xrightarrow{\phi_2}
P^{b_1}\xrightarrow{\phi_1} P \xrightarrow{} A\xrightarrow{}0\, ,\\
0\xrightarrow{}
{\bar{P}}^{b_{e-2}}\xrightarrow{{\bar{\phi}}_{e-2}}\dots
\xrightarrow{} {\bar{P}}^{b_2}\xrightarrow{{\bar{\phi}}_2}
{\bar{P}}^{b_1}\xrightarrow{{\bar{\phi}}_1} {\bar{P}} \xrightarrow{}
{\bar{A}}\xrightarrow{}0\, ,\\
\end{gather*}
so that
\begin{list}{\textup{(\roman{temp})}}{\usecounter{temp}}
        \item the second resolution is the associated graded complex
attached to
                the first;
        \item ${\bar{\phi}}_i$ is homogeneous, of degree $1$ ($i>1$) or $2$
($i=1$);
        \item $b_i=i\binom{e-1}{i+1}$.
\end{list}
\end{theorem}

Actually we will only need part (i) and (ii) for $i=1,2,3$. These
imply that the ring of a rational surface singularity is something we
call a QL ring (see Section~\ref{QL_rings}), making it easy to
compare the equations and relations defining the blow up locally with
those of $X$.

The only result from previous work on $T^2$ we need is about the
module structure. This is \cite[Theorem 5.1.1 (1)]{bc:hypr}, but the
statement in that paper about annihilators of $T^2$ is incorrect.
What actually is proven is
\begin{theorem}[Behnke--Christophersen]\label{BC}
If $X$ is a rational surface singularity with $e\ge 4$ and $x\in
{\mathfrak m}\setminus{\mathfrak m}^2$ is generic (i.e. projects onto
a general element of ${\mathfrak m}/{\mathfrak m}^2$), then $\dim
T^2_X/{\mathfrak m} T^2_X = (e-2)(e-4)$ and $xT^2_X = {\mathfrak m}
T^2_X$.
\end{theorem}

In other words $\dim T^2_X/xT^2_X = (e-2)(e-4)$.  This will be
important in Section~\ref{t2_formula}.

\subsection{Cotangent cohomology}\label{cotangent} We review some
properties of the cotangent complex. For our later use it is enough
to assume that we have a noetherian
ring $S$ and an $S$ algebra $A$ of essentially finite type. There
exists a complex of free $A$ modules; the {\it cotangent complex}
$\mathbb L_{\cdot}^{A/S}$. See \cite[p. 34]{an:hom} for a definition.
For an $A$ module $M$ we get the {\it cotangent cohomology} modules
$T^i(A/S;M):=H^i(\Hom_A(\mathbb L_{\cdot}^{A/S},M))$. If $S$ is the
ground field we abbreviate $T^i(A/S;M)=:T^i_A(M)$ and
$T^i_A(A)=:T^i_A=:T^i_X$ if $X=\Spec A$.

The first three modules are important in deformation theory and we
could have given an ad hoc definition as follows.  Let $P$ be a
polynomial
$S$ algebra (or the localization of such an algebra) mapping onto $A$
so that $A\simeq P/I$ for an ideal $I$. Let
$$0 \to R \to F \stackrel{j}{\to} P \to A \to 0$$
be an exact sequence presenting $A$ as a $P$ module with $F\simeq
P^m$ free. We have $T^0(A/S;M) =\Der_S(A,M)$, the module of $S$
derivations into $M$. The cokernel of the natural map $\Der_S(P,M)
\to \Hom_A(I/{I^2},M)$ is $T^1(A/S;M)$. Let $R_0$ be the sub-module
of $R$ generated by the trivial relations; i.e. those of the form
$j(x)\, y -j(y)\, x$. Then $R/{R_0}$ is an $A$ module and we have an
induced map
$\Hom_A(F/{R_0} \otimes_{P} A,M) \to \Hom_A(R/{R_0},M)$. The cokernel
is $T^2(A/S;M)$. Notice that $\Hom_A(F/{R_0} \otimes_{P} A,M)$ is
just the sum of $m$ copies of $M$ and the map is $$(\alpha_1,\dots
,\alpha_m) \mapsto \lbrack \bar r \mapsto \sum r_i\, \alpha_i
\rbrack$$ where $r\in F$ represents $\bar r \in R/{R_0}$.

If $A$ is a {\it smooth} $S$ algebra then $T^i(A/S;M)=0$ for
$i \geq 1$ and all $A$ modules $M$. As usual a short exact sequence
of $A$ modules induces a long exact sequence in the $T^i(A/S;\star)$.
More importantly, two ring homomorphisms $S \to A \to B$ induce a
long
exact sequence called the {\it Zariski--Jacobi} long exact sequence;
$$\cdots \to T^i(B/A;M) \to T^i(B/S;M) \to T^i(A/S;M) \to
T^{i+1}(B/A;M)  \to \cdots$$
where $M$ is a $B$ module. (See \cite[Th{\'e}or{\`e}me 5.1]{an:hom}.)

If $Y$ is a scheme we may globalise the above local construction.
(See for example \cite[Appendice]{an:hom}, \cite[2.2.3]{bu:the} and
\cite[3.2]{la:for} for details and proofs.) If $\mathcal S$ is a
sheaf
of rings and $\mathcal A$ an $\mathcal S$ algebra we set ${\mathcal
L}_{\cdot}^{{\mathcal A}/{\mathcal S}}$ to be the complex of sheaves
associated with the presheaves $U\mapsto {\mathbb
L}_{\cdot}^{{\mathcal A}(U)/{\mathcal S(U)}}$. Let $\mathcal F$ be an
${\mathcal A}$ module. We get the cotangent cohomology sheaves
${\mathcal T}_{{\mathcal A}/{\mathcal S}}^i({\mathcal F})$ as the
cohomology sheaves of ${\mathcal Hom}_{\mathcal A}({\mathcal
L}_{\cdot}^{{\mathcal A}/{\mathcal S}},{\mathcal F})$ and the
cotangent cohomology groups $T_{{\mathcal A}/{\mathcal
S}}^i({\mathcal F})$ as the cohomology of $\Hom_{\mathcal
A}({\mathcal L}_{\cdot}^{{\mathcal A}/{\mathcal S}},{\mathcal F})$.

Because of the functoriality of these constructions: ${\mathcal
T}_{{\mathcal A}/{\mathcal S}}^i({\mathcal F})$
is the sheaf associated to the presheaf $U\mapsto T^i({\mathcal
A}(U)/{\mathcal S}(U);{\mathcal F}(U))$ and $T_{{\mathcal
A}/{\mathcal S}}^{\star}({\mathcal F})$ is the hyper-cohomology of
${\mathcal Hom}_{\mathcal A}({\mathcal L}_{\cdot}^{{\mathcal
A}/{\mathcal S}},{\mathcal F})$. In particular there is a
``local-global" spectral sequence
\begin{equation}\label{local-global}
H^p(Y,{\mathcal T}_{{\mathcal A}/{\mathcal S}}^q({\mathcal
F}))\Rightarrow T_{{\mathcal A}/{\mathcal S}}^n({\mathcal F})\, .
\end{equation}
If ${\mathcal A}$ is the structure sheaf ${\mathcal O}_Y$ and $S$ is
the ground field, then (abbreviating as
above) the $T^i_Y$ play a role in the deformation theory of $Y$
similar to the local case. The spectral sequence~\eqref{local-global}
becomes $H^p(Y,{\mathcal T}_Y^q)\Rightarrow T_Y^n$ and shows how
local and global deformations contribute to the total.

\subsection{Cotangent cohomology and modifications of rational
singularities} Consider first any morphism of schemes $f:Y\to X$. Let
$f^{-1}$ be the sheaf theoretical adjoint functor of $f_{\star}$ as
defined in for example \cite[II,1]{ha:alg}. We have the for us very
important result in \cite[Appendice. Proposition 56]{an:hom} which we
translate to our notation. (Notice that the $f^{\star}$ in
\cite{an:hom}
is $f^{-1}$ and not the $f^{\star}$ in standard algebraic geometry
notation).
\begin{proposition}\label{if_R1=0} If $f:Y\to X$ is a morphism of
schemes,
${\mathcal A}$ is an ${\mathcal S}$ algebra on $X$ and $\mathcal F$
is an $f^{-1}{\mathcal A}$ module on $Y$ such that
$R^kf_{\star}(\mathcal F)=0$ for $k\ge 1$, then there are natural
isomorphisms
$$T_{f^{-1}{\mathcal A}/f^{-1}{\mathcal S}}^i({\mathcal F})\simeq
T_{{\mathcal A}/{\mathcal S}}^i(f_{\star}{\mathcal F})$$
for all $i\ge 0$.
\end{proposition}

Assume now that $X=\Spec A$ is a normal singularity and $f:Y\to X$ is
a modification; i.e. f is proper and birational. We will slightly
abuse notation and write $f^{-1}A$ for $f^{-1}{\mathcal O}_X$. From
the spectral sequence~\eqref{local-global} and
Proposition~\ref{if_R1=0} we get immediately
\begin{theorem}\label{spectral}
If $f:Y\to X=\Spec A$ is a modification and $\mathcal F$ is a
coherent sheaf on $Y$ with $R^kf_{\star}(\mathcal F)=0$ for $k\ge 1$,
then there is a spectral sequence $\{E_r^{p,q}\}$ with
$E_2^{p,q}=H^p(Y,{\mathcal T}_{f^{-1}A}^q({\mathcal F}))$ such that
$$E_2^{p,q}\Rightarrow T_A^n(f_{\star}{\mathcal F})\, .$$
\end{theorem}

\begin{remark} The results we have compiled from the literature to
get Theorem~\ref{spectral} involve injective resolutions to compute
hyper-cohomology etc.. In the computational part of this paper it
will be important to know some of the maps from the spectral sequence
explicitly, and therefore in terms of \v{C}ech cohomology. We will
state these descriptions without proof. For a proof of
Theorem~\ref{spectral} using \v{C}ech cohomology (done before we
found the relevant known results) and explaining the maps see
\cite{gu:alg}.
\end{remark}

\begin{corollary}\label{ex-sequence}
If $f:Y\to X=\Spec A$ is a modification of a rational surface
singularity then there are exact sequences
$$0\to H^1( Y, {\mathcal T}_{f^{-1}A}^{i-1}({\mathcal O}_Y))\to T^i_X
\to H^0( Y, {\mathcal T}_{f^{-1}A}^i({\mathcal O}_Y))\to 0 $$
for all $i\ge 0$.
\end{corollary}
\begin{proof} The spectral sequence in Theorem~\ref{spectral} is
derived from a double complex for computing the hyper-cohomology
$T_{f^{-1}A}^n ({\mathcal O}_Y)$. On the other hand, since $X$ is
affine, $H^i(Y,{\mathcal G})=0$ for $i\ge 2$ for any coherent
${\mathcal G}$. Thus the $E_2^{p,q}$ consists of two adjacent
non-zero rows. The result follows from standard arguments. Notice
that $f_{\star}{\mathcal O}_Y\simeq {\mathcal O}_X$ by normality.
\end{proof}

We will use Corollary~\ref{ex-sequence}  when the modification is the
blow-up $\pi:\widehat{X}\to X$ to prove our formulas. To shorten
notation
we define the sheaves on $\widehat{X}$
$${\mathcal F}^i:={\mathcal T}_{\pi^{-1}A}^i({\mathcal
O}_{\widehat{X}})\, .$$
Notice that there are natural maps ${\mathcal T}_Y^i\to {\mathcal
T}_{f^{-1}A}^i({\mathcal O}_Y)$ that induce, with the help of the
exact sequences, natural maps $T^i_Y\to T^i_X$. If $i=1$ these are
the tangent maps to the contraction of deformations of $Y$ to
deformations of $X$. They behave very sporadically, and we have not
found them useful for proving the formulas. Instead we will have to
make some unnatural maps relating the ${\mathcal F}^i$ to the
$T^i_{\widehat{X}}$.

\section{QL-rings and blowing up}

\subsection{Associated graded rings and standard bases} We recall
some facts
regarding associated graded rings and standard bases.
\begin{definition} Let $R$ be a noetherian local ring with maximal
ideal $\mathfrak m$
and let $M$ be a finitely generated $R$-module. Let $N\subset M$ be a
submodule. We set
\begin{align*}
G(\mathfrak m, M) &:= \underset{i\ge 0}{\oplus}\mathfrak m^iM/
        \mathfrak m^{i+1}M\\
G(\mathfrak m, N\subset M) &:= \underset{i\ge 0}{\oplus}\mathfrak
m^iM\cap N +
        \mathfrak m^{i+1}M/\mathfrak m^{i+1}M\\
B(\mathfrak m, M) &:= \underset{i\ge 0}{\oplus}\mathfrak m^iM\, .
\end{align*}
Also, for any nonzero $m\in M$ we put
$$ \ord(\mathfrak m,M)(m)=\sup\left\{n|m\in \mathfrak m^n
M\right\}, $$
and if $m\not = 0$ and $\ord(\mathfrak m,M)(m)=d$ we define the
initial form
$$\initform(\mathfrak m, M)(m) = m + \mathfrak m^dM/\mathfrak
m^{d+1}M $$
Let $m_1,\dots,m_t$ be generators for $N$. Then $m_1,\dots,m_t$ is
called a
{\it standard basis} for $N$ if the submodule $G(\mathfrak m,
N\subset M)$ is
generated by
$\initform(\mathfrak m, M)(m_1),\dots,\initform(\mathfrak m, M)(m_t).$
\end{definition}
We will write $\ord(m)$ for $\ord(\mathfrak m,M)(m)$ and
$\initform(m)$ for $\initform(\mathfrak m, M)(m)$ when no
misunderstanding is
likely to occur.

We will need the following result -- see e.g. \cite[Theorem
13.7]{hio:equ}.
\begin{theorem} \label{about_standard_base}
Let $R$ be a noetherian local ring
with maximal ideal $\mathfrak m,$ let $M$ be a finitely generated
$R$-module and
let $N\subset M$ be a submodule. Then $m_1,\dots, m_t$ is a
standard basis for $N$
if and only if for any $z\in N$ there are $a_1,\dots,a_t\in R$ such
that
$z=a_1m_1+\cdots a_tm_t$ and $\ord(z)\le \ord(a_i)+\ord(m)$ for all
$i.$
\end{theorem}

\subsection{QL-rings} \label{QL_rings}
Let $P=\mathbb C\left[x_1,\dots,x_e\right]_{\mathfrak m}$ be the
polynomial ring with $e$ generators localized in $\mathfrak m =
(x_1,\dots,x_e)$. We let $\mathfrak m$ denote both the maximal ideal
in $P$ and quotients of $P$ unless this causes confusion.
\setcounter{temp}{0}
\begin{definition}
\label{def_of_QL}
 We will say that $A$ is a {\it QL-ring} (quadratic generators and
linear relations) if $A=P/I$ where $I\subset P$ is a prime ideal such
that
\begin{list}{\textup{(\roman{temp})}}{\usecounter{temp}}
\item  The ideal $I$ has a standard basis $f_1,\dots,f_m$ with
        $\ord(f_i)=2$ for $i=1,\dots ,m$.
\item The relation module
        $R=\left\{(p_1,\dots,p_l)\in P^m|\sum p_if_i=0\right\}$ has a
standard
        basis $r^1,\dots,r^s$ with $\ord(r^i)=1$ for
        $i=1,\dots ,s$.
\item The $\initform(f_i)$ and $\initform(r^i)$ are linearly
independent.
\end{list}
We say that $X=\Spec A$ is a {\it QL singularity} if $A$ is a QL-ring.
\end{definition}

Rational surface singularities with $e\ge 4$ are QL singularities by
Theorem~\ref{Wah}. Another example is the class of minimal elliptic
surface singularities with $e\ge 5$ (\cite[Theorem 2.8]{wa:equ}).

QL-singularities have an algebraic property that will be very
important for us in the proof of the formulas. We state it here for
future reference.
\setcounter{temp}{0}
\begin{lemma}\label{ql_lemma}
Suppose $A$ is a QL-ring and assume $f_1,\dots,f_m$ and
$r^1,\dots,r^s$ are standard bases as in Definition~\ref{def_of_QL}.
Let $x\in {\mathfrak m} \setminus {\mathfrak m}^2$.
\begin{list}{\textup{(\roman{temp})}}{\usecounter{temp}}
\item Every $f_j$ is involved in some relation $r=(r_1,\dots ,r_m)$
with
        $r_j\in {\mathfrak m}\setminus {\mathfrak m}^2$.
\item Every $f_j$ is involved in some relation $r=(r_1,\dots ,r_m)$
with
        $r_j\notin (x)$.
\end{list}
\end{lemma}
\begin{proof} The first statement is proven in \cite[2.5]{wa:jac} and
the second one follows from the same argument, so we repeat it here.
We may assume $j=1$ and consider the trivial relation
$(f_2,-f_1,0,\dots ,0)=\sum q_ir^i$. Thus all $q_i\in{\mathfrak m}$
and $f_2=\sum q_ir^i_1\notin {\mathfrak m}^3$. Also -- $f_2\notin
(x)$ since $x\notin {\mathfrak m}^2$ and $A$ is a domain. So some
$r^i_1\notin {\mathfrak m}^2$ and some $r^j_1\notin (x)$.
\end{proof}

As a consequence we get a slight generalization of \cite[Lemma
2.6]{wa:jac}.
\begin{lemma}\label{hom_in_m}
If $A=P/I$ is a QL-ring and $x\in ({\mathfrak m} \setminus {\mathfrak
m}^2)\cup \{0\}$, then every $\phi\in\Hom_A(I/I^2,A/(x))$ has $\Image
(\phi)\subseteq {\mathfrak m}_{A/(x)}$.
\end{lemma}
\begin{proof} Assume $\phi(f_1)=1$ and let $h_i\in P$ represent
$\phi(f_i)\in A/(x)$. After changing $f_i$ to $f_i-h_if_1$ we may
assume $\phi(f_i)=0$ for $i\ge 2$. (Remember that the
$\initform(f_i)$ are independent.) This contradicts
Lemma~\ref{ql_lemma}, since for all relations $\sum
r_i\phi(f_i)=r_1\equiv 0\mod (x)$.
\end{proof}

\subsection{QL singularities and blowing up}
Let $A=P/I$ be a QL-ring. Let $\pi:\widehat{X}=\Proj B(\mathfrak
m,A)\to \Spec A$ be the blow up of $\Spec A$. Recall that the
blow-up $\widehat{X}$ is covered by $\Spec B(\mathfrak m,A)_{(xt)}$
for
$x\in\mathfrak m_A\subset A,$ where we consider $B(\mathfrak m,A)$ as $A[\mathfrak mt].$

For $x\in P,x\notin I$  we denote by $P\left[\frac{\mathfrak
m}{x}\right]$ the subring of $\mathbb C(x_1,\dots,x_e)$ generated
by the image $P$ and the elements of the form $\frac{a}{x}$ where
$a\in
\mathfrak m$. Similarly we denote by $A\left[\frac{\mathfrak
m}{x}\right]$ the subring of $A$'s quotient field generated by the
image of
$A$ and the elements of the form $\frac{a}{x}$ where $a\in \mathfrak
m$.
For the covering of the blow-up we have $B(\mathfrak
m,A)_{(xt)}\simeq A\left[\frac{\mathfrak m}{x}\right]$
and the restriction of $\pi$ to $\Spec B(\mathfrak m,A)_{(xt)}$ is
induced by the inclusion $A\subset A\left[\frac{\mathfrak
m}{x}\right]$.

\setcounter{temp}{0}
\begin{lemma}\label{IR_lemma}
Suppose $A=P/I$ is a QL-ring and assume $f_1,\dots,f_m$ and
$r^1,\dots,r^s$ are standard bases as in Definition~\ref{def_of_QL}.
Then
\begin{list}{\textup{(\roman{temp})}}{\usecounter{temp}}
\item $A\left[\frac{\mathfrak m}{x}\right]\simeq
        P\left[\frac{\mathfrak m}{x}\right]/I_B$ where $I_B$ is generated
by
        ${f_1}/{x^2},\dots,{f_m}/{x^2}$.
\item  The relation module
        $$S=\left\{(p_1,\dots,p_m)\in P\left[\frac{\mathfrak
m}{x}\right]^m|
        \sum p_i\frac{f_i}{x^2}=0 \;\operatorname{in}
        \; P\left[\frac{\mathfrak m}{x}\right]\right\}$$
        is generated by ${r^1}/{x},\dots,{r^s}/{x}$.
\end{list}
\end{lemma}
\begin{proof}
The first statement is a special case of a well known property of
blow-ups, see e.g. \cite[Proposition 13.13]{hio:equ}.

We prove the second statement for lack of reference. Let $R$ be the
relation module for the $f_i$, and set $\frac{R}{x}$ to be the
$P\left[\frac{\mathfrak m}{x}\right]$ module generated by
${r^1}/{x},\dots,{r^s}/{x}$. Clearly $\frac{R}{x}\subset S$. Choose
some $p=(p_1,\dots,p_m)\in S$. We may find an $n$ such that
$x^{n-2}p_i\in P$ for all $i=1,\dots,m$. Thus $\sum x^np_if_i=0$, so
$x^np\in R$. Since $x^np\in R$, we are by
Theorem~\ref{about_standard_base} able to find $q_1,\dots,q_s$ such
that
$x^np=\sum q_jr^j$ and $\ord(\mathfrak m, P)(q_j)+\ord(\mathfrak
m, P^m)(r^j)\ge\ord(\mathfrak m, P^m)(x^np)$. Now $\ord(r^j)=1$
and $\ord(x^np)\ge n$, hence $\ord(q_j) \ge n-1$.  We end up with $p
= \sum ({q_j}/{x^{n-1}})(r^j/{x})$ with ${q_j}/{x^{n-1}}\in
P\left[\frac{\mathfrak m}{x}\right]$, which shows that
$p\in\frac{R}{x}$.
\end{proof}

\section{The formulas}

\subsection{Computation of the ${\mathcal F}^i$.} \label{Computation}

Let $A=P/I$ be a QL-ring, where $P$ is as above with $e=\dim
{\mathfrak m}_A/{\mathfrak m}_A^2$. Let
$\widehat{X}\stackrel{\pi}{\rightarrow} X$ be the blow up of $X=\Spec
A$
and $M$ the blowup of $\Spec P$, so that we may view $\widehat{X}$ as
the
strict transform of $X$ in $M$. Let $C\subset \widehat{X}$ be the
exceptional divisor defined by ${\mathfrak m}{\mathcal
O}_{\widehat{X}}$.
Thus $C=\Proj G(\mathfrak m, A)$.

\begin{nota}\label{standardnot}
Throughout the following proofs we will be working locally on
$\widehat{X}$ with affine charts $\Spec B$ with
$B=A\left[\frac{\mathfrak m}{x}\right]$ as in Lemma~\ref{IR_lemma}.
Set $P_B=P\left[\frac{\mathfrak m}{x}\right]$. We use variables
$x,x_1,\dots ,x_n$ ($n=e-1$) for $P$, so $t_i:={x_i}/{x}$ and $x$
generate $P_B$. Generators for $I$ will be denoted $f_1,\dots ,f_m$
and $g_1,\dots ,g_m$ will be generators of $I_B$ as in
Lemma~\ref{IR_lemma}; i.e. $x^2 g_i= f_i(x,t_1,\dots ,t_n)$ . We view
the $x_i=xt_i$ as elements of $P_B$ as well.
\end{nota}

Let ${\mathcal N}_{{\widehat{X}}/M}$ be the normal sheaf of
$\widehat{X}$ in
$M$ and set $N_A({\mathcal O}_{\widehat{X}}) = {\mathcal
T}^1_{\pi^{-1}A/\pi^{-1}P}({\mathcal O}_{\widehat{X}})$. Let
$\Theta_{\widehat{X}}$ be the tangent sheaf on $\widehat X$. Let
$\Der_C(\widehat{X})$ be
the subsheaf of $\Theta_{\widehat{X}}$ consisting of derivations $D$
with
$D({\mathcal I}_C)\subseteq {\mathcal I}_C$. Finally define
${\mathcal A}^1_{C/{\widehat{X}}}$
to be the cokernel of the map $\Theta_{\widehat{X}}\rightarrow
{\mathcal
O}_C(C)$ defined locally -- where $C$ is defined by $x$ -- as
$D\mapsto D(x)\otimes\frac1x\mod (x)$. Notice that there is an exact
sequence
\begin{equation}\label{A-exseq}
0\rightarrow \Der_C(\widehat{X}) \rightarrow \Theta_{\widehat{X}}
\rightarrow {\mathcal O}_C(C) \rightarrow {\mathcal
A}^1_{C/{\widehat{X}}} \rightarrow 0
\end{equation}
with the maps as above.

First we compare the ${\mathcal F}^i$ with the ${\mathcal
T}^i_{\widehat{X}}$.
\setcounter{temp}{0}
\begin{proposition}
\label{identities}
If $X$ is a QL-singularity, then
\begin{list}{\textup{(\roman{temp})}}{\usecounter{temp}}
        \item ${\mathcal F}^0 \simeq \Der_C(\widehat{X})(C)$.
        \item $N_A({\mathcal O}_{\widehat{X}}) \simeq {\mathcal
        N}_{{\widehat{X}}/M}(2C)$.
        \item There is an exact sequence $0\rightarrow {\mathcal
        A}^1_{C/{\widehat{X}}}(C)\rightarrow {\mathfrak m}{\mathcal F}^1
\rightarrow {\mathcal
        T}^1_{\widehat{X}}(C) \rightarrow 0$.
        \item ${\mathfrak m}{\mathcal F}^2 \simeq {\mathcal
T}^2_{\widehat{X}}$.
\end{list}
The isomorphisms and maps are non-canonical.
\end{proposition}

\begin{proof}
Consider an affine chart $\Spec B$ of $\widehat{X}$ with
$B=A\left[\frac{\mathfrak m}{x}\right]$ as above.
The isomorphism in (i) is given by the map
$\Der_C(\widehat{X})(C)\rightarrow {\mathcal F}^0$ which is locally
$D\otimes\tfrac1x\mapsto\tfrac1x D\rvert_A$. This is easily checked
to be injective, and if $\delta \in \Der(A,B)$ then it comes from a
derivation $D$ determined by $D(x)=x\delta(x)$ and
$D(t_i)=\delta(x_i)-t_i\delta(x)$.

Let $\mathcal J$ be the ideal sheaf of $\widehat X$ in $M$ and $V$
the exceptional divisor of $\pi:M\to \Spec P$. Lemma~\ref{IR_lemma}
implies that $\pi^\star I\simeq {\mathcal J}(-2V)$; this induces the
isomorphism in (ii).

The chain rule and the fact that $x$ is not a zero divisor, yields
the following equalities in $B$:
\begin{subequations}\label{chainrule}
\begin{align}
        \frac{\partial f_j}{\partial x}& = x(x\frac{\partial g_j}{\partial
x}-\sum_i
        t_i\frac{\partial g_j}{\partial t_i})\\
        \frac{\partial f_j}{\partial x_i}& = x\frac{\partial g_j}{\partial
t_i}\, .
\end{align}
\end{subequations}

The isomorphism (ii) takes $\phi\in\Hom_P(I,B)$ to the morphism
determined by $g_j\mapsto \phi(f_j)$. In particular \eqref{chainrule}
shows that it induces (after a twist) a map ${\mathcal
F}^1(-C)\rightarrow {\mathcal T}^1_{\widehat{X}}(C)$ which must be
surjective. We claim that this factors through ${\mathfrak
m}{\mathcal F}^1$. Indeed, if $x[\phi]=0$ in $T^1_A(B)$, then
$x\phi(f_j)=b\,{\partial f_j}/{\partial x}+\sum b_i\,{\partial
f_j}/{\partial x_i}$, so by \eqref{chainrule},
$\phi(f_j)=b\,(x\,{\partial g_j}/{\partial x}-\sum t_i\,{\partial
g_j}/{\partial t_i})+\sum b_i\,{\partial g_j}/{\partial t_i}$. Thus
$[\phi]\otimes x\mapsto 0$. This gives the right surjection in (iii).

Let $K$ be the kernel of this map locally, i.e. of
$xT^1_A(B)\rightarrow T^1_B$. We have $x[\phi]\in K$ iff
$\phi(f_j)=b\,{\partial g_j}/{\partial x}+\sum b_i\,{\partial
g_j}/{\partial t_i}$. But then $x\phi(f_j)=bx\,{\partial
g_j}/{\partial x}+\sum b_ix\,{\partial g_j}/{\partial
t_i}=xb\,{\partial g_j}/{\partial x}+\sum b_i\,{\partial
f_j}/{\partial x_i}$, so $x[\phi]$ equals the class of the map
$f_j\mapsto bx\,{\partial g_j}/{\partial x}$ in $T^1_A(B)$. In
particular $K$ is a cyclic $B$ module generated by the class of the
map $f_j\mapsto x\,{\partial g_j}/{\partial x}$.

This yields a surjection $B\twoheadrightarrow K$. The kernel of this
map is
\begin{multline*}
\{b\in B : bx\,\frac{\partial g_j}{\partial x}=b_0\,\frac{\partial
f_j}{\partial x}+\sum b_i\,\frac{\partial f_j}{\partial x_i}\text{
for
some } b_i\in B,\,j=1,\dots ,m\}\\
=\{b\in B : b\,\frac{\partial g_j}{\partial
x}=b_0x\,\frac{\partial g_j}{\partial x}+\sum
(b_i-b_0t_i)\,\frac{\partial
g_j}{\partial t_i}\}={\mathfrak a} + (x)\, .
\end{multline*}
 where ${\mathfrak a}=\{b : b\,{\partial g_j}/{\partial x}=
 \sum b_i\,{\partial g_j}/{\partial t_i}\}$.
 But clearly this last condition is the same as $b=D(x)$ for
some $D\in \Der(B)$. This gives an exact sequence
$$0\rightarrow B/{{\mathfrak a}+(x)} \rightarrow xT^1_A(B)\rightarrow
T^1_B \rightarrow 0$$
which globalizes to the one in (iii).

Let $R$ and $S$ be as in Lemma~\ref{IR_lemma} and let $R_0$ and $S_0$
be the submodules of Koszul relations. Thus $R\otimes_P P_B \simeq
x\cdot S$ and $R_0\otimes P_B \simeq x^2\cdot S_0$. Now $\Hom_A
(R/R_0,B)$ is the kernel of the natural map $\Hom_A
(R/IR,B)\rightarrow \Hom_A (R_0/IR,B)$, so it is isomorphic to the
kernel of $\Hom_{P_B} (x\cdot S,B)\rightarrow \Hom_{P_B} (x^2\cdot
S_0,B)$. This kernel is again isomorphic to $\Hom_B (S/S_0,B)$ since
$S_0/xS_0$ is annihilated by the non-zero divisor $x$.

This isomorphism induces a surjection $T^2_A(B) \twoheadrightarrow
T^2_B$. One checks that the kernel is $\{ [\phi] \in T^2_A(B) :
x\cdot [\phi] \equiv 0\}$ which is also the kernel of the
multiplication map $T^2_A(B) \stackrel{\cdot x}{\rightarrow}
T^2_A(B)$. This induces the isomorphism (iv) locally.
\end{proof}

\subsection{The $T^2$ formula.} \label{t2_formula}

\begin{proposition}
\label{h0-zero}
If $X$ is a rational surface singularity, then
$H^0(\widehat{X},{\mathcal
F}^1\rvert_C) = H^0(\widehat{X},{\mathcal F}^2\rvert_C)=0$.
\end{proposition}
\begin{proof} From the quotient map ${\mathcal
O}_{\widehat{X}}\rightarrow {\mathcal O}_C$ and
Theorem~\ref{spectral} we
get the following commutative diagram with surjective horizontal maps;
$$\begin{CD}
T^2_X @>>> H^0(\widehat{X},{\mathcal F}^2) \\
@V\alpha VV @VV\beta V\\
T^2_X(\pi_{\star}{\mathcal O}_C) @>>>H^0(\widehat{X},{\mathcal
T}^2_{\pi^{-1}A}({\mathcal O}_C))\\
\end{CD}$$
Now $\pi_{\star}{\mathcal O}_C\simeq A/{\mathfrak m}$. For a rational
singularity with $e\ge 5$ the ``relations among relations" are
generated by independent linear ones (Theorem~\ref{Wah}). We may
argue as in Lemma~\ref{hom_in_m} to show that the images of all
$\phi\in\Hom_A(R/R_0,A)$ are in $\mathfrak m$. So $\alpha$ is the
zero-map and therefore $\beta$ is the zero-map. On the other hand
$\beta$ factors $H^0({\mathcal F}^2)\rightarrow H^0({\mathcal
F}^2\rvert_C) \rightarrow H^0({\mathcal T}^2_{\pi^{-1}A} ({\mathcal
O}_C))$. The second map is injective and the cokernel of the first
map is contained in $H^1({\mathfrak m}{\mathcal F}^2)$ which is zero
by Proposition~\ref{identities}. This proves that
$H^0(\widehat{X},{\mathcal F}^2\rvert_C)=0$.

In the case of ${\mathcal F}^1\rvert_C$ we can make a direct
calculation relying only on the QL property. As above we consider the
injective map $H^0({\mathcal F}^1\rvert_C) \hookrightarrow
H^0({\mathcal T}^1_{\pi^{-1}A} ({\mathcal O}_C))$.
Since the ${\pi^{-1}A}$ module structure on ${\mathcal O}_C$ is
defined by $A\twoheadrightarrow A/{\mathfrak m}\simeq {\mathbb
C}\hookrightarrow {\mathcal O}_C$, we have  ${\mathcal
T}^1_{\pi^{-1}A} ({\mathcal
O}_C)\simeq m{\mathcal O}_C$ where $m$ is the minimal number of
generators for $I$. In particular $H^0({\mathcal T}^1_{\pi^{-1}A}
({\mathcal O}_C))\simeq {\mathbb C}^m$. A global section of
$H^0({\mathcal F}^1\rvert_C)$ must therefore be locally represented
by a homomorphism that looks like $f_j\mapsto \lambda_j + {\mathcal
I}_C$ with
$\lambda_j\in\mathbb C$.

We claim that for every $f_j$ there exists a chart with coordinate
ring $B = A\left[\frac{\mathfrak m}{x}\right]$, such that there are
no $\phi\in\Hom_A(I/I^2,B)$ with $\phi(f_j)\equiv\lambda \mod (x)$
and $\lambda\ne 0$ a constant. To prove this consider for $f_j$ a
relation as in Lemma~\ref{ql_lemma} and set $x=r_j$. If any of the
other $r_k\in (x)$, say $r_k=h_kx$, change $f_j$ to $f_j+\sum_k
h_kf_k$. Thus we may assume all other $r_k\in {\mathfrak m}\setminus
(x)$. We must have $\sum r_i\phi(f_i) = 0$ in $B$, so
$\phi(f_j)=-\sum_{i\ne j}(r_i/x)\phi(f_i)$ in $B$, but by the
assumption on these $r_i$, none of the $(r_i/x)$ are constants.
\end{proof}

The following result follows immediately from
Proposition~\ref{identities} and Proposition~\ref{h0-zero}.
\begin{corollary}\label{f2t2}
If $X$ is a rational surface singularity then
$$H^0(\widehat{X},{\mathcal F}^2) \simeq H^0(\widehat{X},{\mathcal
T}^2_{\widehat{X}})$$
and the sequence
$$0\rightarrow H^1(\widehat{X},{\mathfrak m}{\mathcal F}^1)
\rightarrow
H^1(\widehat{X},{\mathcal F}^1) \rightarrow H^1(\widehat{X},{\mathcal
F}^1\rvert_C) \rightarrow 0$$
is exact.
\end{corollary}

Let us now concentrate on $H^1(\widehat{X},{\mathcal F}^1)$. Using
Theorem~\ref{BC} we will prove via two lemmas that
$\dim_{\mathbb C} H^1(\widehat{X},{\mathcal
F}^1\rvert_C)=(e-2)(e-4)$. (We view cohomology groups on
$\widehat{X}$ as $A$ modules by their isomorphisms with
$R^i\pi_{\star}$'s).

\begin{lemma}\label{socledim}
Suppose $X$ is a rational surface singularity and $x\in\mathfrak m$ is
generic.  If $K$ is a submodule of $T^2_X$ containing the kernel of
multiplication by $x$, then $\dim_{\mathbb C} K/xK = (e-2)(e-4)$.
\end{lemma}
\begin{proof}
Let $L$ be the kernel of the multiplication map $K\stackrel{\cdot
x}{\rightarrow} K$ and $M$ the kernel of multiplication by $x$ on $T^2_X$.
Clearly $L=M\cap K$, but we have assumed that $M\subseteq K$, so
$L=M$.  Thus $\dim_{\mathbb C} K/xK = \dim_{\mathbb C} T^2/xT^2$
which equals $(e-2)(e-4)$ by Theorem~\ref{BC}.
\end{proof}

\begin{lemma}\label{socle_in_h1}
Suppose $X$ is a rational surface singularity and $x\in\mathfrak m$ is
generic.  If we identify $H^1(\widehat{X},{\mathcal F}^1)$ with the
kernel of $T^2_X\rightarrow H^0(\widehat{X},{\mathcal F}^2)$, then
$H^1(\widehat{X},{\mathcal F}^1)$ contains the kernel of
multiplication by $x$.
\end{lemma}
\begin{proof}
We may assume by genericity of $x$ that a global section vanishes if
and only if it vanishes in the chart $\Spec A\left[\frac{\mathfrak m}{x}\right]$.  (We know for example
that $H^0({\mathcal F}^2)\simeq H^0({\mathcal T}^2)$ and ${\mathcal
T}^2$ has support at points.)  Thus we need to show that the local
maps $T^2_A\rightarrow T^2_A(B)$ map an element of the kernel of
multiplication by $x$ to zero if $B=A\left[\frac{\mathfrak
m}{x}\right]$.  If $[\phi]$ is in this kernel we may find $a_1,\dots
,a_m\in A$ such that $x\phi(r)=\sum a_ir_i$ in $A$ for all relations
$r=(r_1,\dots ,r_m)$.  We claim that $a_1,\dots ,a_m\in \mathfrak m$.
This is because the application $f_j\mapsto a_j$ defines an element of
$\Hom_A\left(I/I^2,A/(x)\right)$, so Lemma~\ref{hom_in_m} applies.
But then $a_i/x\in B$ and $\phi(r)=\sum (a_i/x)r_i$ in $B$, thus
$[\phi]\mapsto 0$.
\end{proof}

\begin{proposition}\label{h1fC}
If $X$ is a rational surface singularity, then $\dim_{\mathbb C}
H^1(\widehat{X},{\mathcal F}^1\rvert_C)=(e-2)(e-4)$.
\end{proposition}
\begin{proof}  Since $x$ is generic, the cokernel of ${\mathcal
F}^1\xrightarrow{\cdot x} {\mathfrak m}{\mathcal F}^1$ has support at
points, so
$xH^1({\mathcal F}^1)\simeq H^1({\mathfrak m}{\mathcal F}^1)$. The
result now follows from Corollary~\ref{f2t2}, Lemma~\ref{socledim}
and Lemma~\ref{socle_in_h1}.
\end{proof}

Putting all of this together we get the formula for $\dim_{\mathbb
C}T^2_X$. We first define the ``correction term". We will see several
other ways of defining this number in Section~\ref{Cor_term}.
\begin{definition}\label{cor}
 If $X$ is a rational surface singularity, we define the invariant
$$c(X):=\dim_{\mathbb C}H^1(\widehat{X},{\mathfrak m}{\mathcal
F}^1)\,.$$
\end{definition}

\begin{theorem}\label{t2}
If $X$ is a rational surface singularity of embedding dimension $e$
and $\widehat{X}$ is the blow up of $X$, then
$$\dim_{\mathbb C}T^2_{X}=(e-2)(e-4) +
\sum_{p\in{\widehat{X}}}\dim_{\mathbb C}T^2_{(\widehat{X},p)} +
c(X)\, .$$
\end{theorem}
\begin{proof} This follows from Corollary~\ref{ex-sequence},
Corollary~\ref{f2t2} and Proposition~\ref{h1fC}.
\end{proof}

\subsection{The $T^1$ formula.} \label{t1_formula}

First notice that Proposition~\ref{identities} implies that the exact
sequence \eqref{A-exseq} extends (after a twist) to an exact sequence
\begin{equation}\label{long-exseq}
0\rightarrow {\mathcal F}^0 \rightarrow \Theta_{\widehat{X}}(C)
\rightarrow {\mathcal O}_C(2C) \rightarrow {\mathfrak m}{\mathcal
F}^1 \rightarrow {\mathcal T}^1_{\widehat{X}}(C) \rightarrow 0\, .
\end{equation}

Only the two first sheaves have support outside $C$; i.e. have an
infinite dimensional $H^0$. On the other hand, the sequence induces
an exact sequence
$$0\rightarrow H^0(\widehat{X},{\mathcal F}^0) \rightarrow
H^0(\widehat{X},\Theta_{\widehat{X}}(C))
\rightarrow H^0(\widehat{X},{\mathcal O}_C(2C))\, .$$
Now $C$ is an arithmetically Cohen-Macaulay curve in ${\mathbb
P}^{e-1}$,(see e.g. \cite{wa:equ}). In particular
$H^0(\widehat{X},{\mathcal O}_C(2C))=H^0(C,{\mathcal O}_C(-2))=0$.
Thus
the sequence~\eqref{long-exseq} induces an isomorphism
$H^0(\widehat{X},{\mathcal F}^0) \simeq
H^0(\widehat{X},\Theta_{\widehat{X}}(C))$. (They are actually
isomorphic to
$\Der(A)$, which is seen immediately from Theorem~\ref{spectral}.)

Proposition~\ref{h0-zero} tells us that $h^0({\mathcal
F}^1)=h^0({\mathfrak m}{\mathcal F}^1)$. Since ${\mathcal
T}^1_{\widehat{X}}$ is a skyscraper sheaf $h^0({\mathcal
T}^1_{\widehat{X}})=h^0({\mathcal T}^1_{\widehat{X}}(C))$ and
$h^1({\mathcal
T}^1_{\widehat{X}}(C))=0$.
Using all this information and Corollary~\ref{ex-sequence} the
sequence~\eqref{long-exseq} yields the formula
\begin{equation}\label{first_h_formula}
\begin{split}
\dim_{\mathbb C}T^1_X& =h^1({\mathcal F}^0)+h^0({\mathcal F}^1) \\
 & = h^1(\Theta_{\widehat{X}}(C)) + \chi({\mathcal O}_C(2C)) +
h^0({\mathcal T}^1_{\widehat{X}}) + c(X)\, .
\end{split}
\end{equation}

Consider now the minimal resolution $\widetilde{X}$ of $X$ which
factors
$\widetilde{X}\stackrel{p}{\to}\widehat{X}\stackrel{\pi}{\to}X$.
Clearly ${\mathcal O}_{\widetilde{X}}(-Z)\simeq
p^{\star} {\mathcal O}_{\widehat{X}}(-C)$, so we also have ${\mathcal
O}_{\widetilde{X}}(kZ)\simeq p^{\star} {\mathcal
O}_{\widehat{X}}(kC)$.

If we use the projection formula on $\Theta_{\widetilde{X}}(Z)$ we
find
that
\begin{equation}\label{proj}
R^ip_{\star}\Theta_{\widetilde{X}}(Z)\simeq
R^ip_{\star}\Theta_{\widetilde{X}}\otimes{\mathcal
O}_{\widehat{X}}(C)\,.
\end{equation}
It is also true for rational surface singularities that
$p_{\star}\Theta_{\widetilde{X}}\simeq \Theta_{\widehat{X}}$
(\cite{bw:loc}), so
in particular, $p_{\star}\Theta_{\widetilde{X}}(Z)\simeq
\Theta_{\widehat{X}}(C)$. Since $H^2$'s vanish, the Leray spectral
sequence gives in our situation, an exact sequence
$$0\rightarrow H^1(\widehat{X},\Theta_{\widehat{X}}(C)) \rightarrow
H^1(\widetilde{X},\Theta_{\widetilde{X}}(Z))
\rightarrow
H^0(\widehat{X},R^1p_{\star}\Theta_{\widetilde{X}}(Z))\rightarrow
0\, .$$
Also by \eqref{proj} we see that $h^0(
R^1p_{\star}\Theta_{\widetilde{X}}(Z)) =
h^0( R^1p_{\star}\Theta_{\widetilde{X}})$.

Consider the exact sequence
\begin{equation}\label{thetaZ}
0\rightarrow \Theta_{\widetilde{X}} \rightarrow
\Theta_{\widetilde{X}}(Z)
\rightarrow \Theta_{\widetilde{X}}\otimes {\mathcal O}_Z(Z)
\rightarrow 0
\, .
\end{equation}
We state and prove for lack of reference the following
\begin{lemma}\label{h0_thetaZ}
If $X$ is a rational surface singularity, then the induced map
$$H^0(\widetilde{X},\Theta_{\widetilde{X}}) \rightarrow
H^0(\widetilde{X},\Theta_{\widetilde{X}}(Z))$$ is an isomorphism.
\end{lemma}
\begin{proof}
There is a well known exact sequence on the resolution of a normal
singularity
$$0\rightarrow \Der_E(\widetilde{X}) \rightarrow
\Theta_{\widetilde{X}}
\rightarrow \bigoplus {\mathcal O}_{E_i}(E_i) \rightarrow 0 $$
(see \cite[Proposition 2.2]{wa:equi}). After tensoring this sequence
with
${\mathcal O}_{\widetilde{X}}(Z)$ and applying $H^0$ we get a
commutative
diagram
$$\begin{CD}
H^0(\widetilde{X},\Der_E(\widetilde{X})) @>>>
H^0(\widetilde{X},\Theta_{\widetilde{X}}) \\
@V\alpha VV @VV\beta V\\
H^0(\widetilde{X},\Der_E(\widetilde{X})(Z)) @>>>
H^0(\widetilde{X},\Theta_{\widetilde{X}}(Z))\\
\end{CD}$$
where all the maps are injective. The sheaves ${\mathcal
O}_{E_i}(E_i)$ and ${\mathcal O}_{E_i}(E_i+Z)$ on $E_i\simeq {\mathbb
P}^1$ have negative degree, so the horizontal maps are also
surjective. The cokernel of $\alpha$ sits in
$H^0(\widetilde{X},\Der_E(\widetilde{X})\otimes {\mathcal O}_Z(Z))$
which is
trivial by a vanishing result -- $H^1_E(\Der_E(\widetilde{X}))=0$ --
of
Wahl. See \cite[Corollary 2.6]{bk:inf} for an argument. So $\alpha$,
and
therefore $\beta$, is an isomorphism.
\end{proof}

\begin{remark}
There is something to prove, since
$H^0(\widetilde{X},\Theta_{\widetilde{X}}\otimes {\mathcal O}_Z(Z))$
is in
general non-trivial. In fact, if $Z$ is reduced, then it has
dimension equal to $\dim_{\mathbb C} H^1_E(\Theta_{\widetilde{X}})$
which
again equals the number of $-2$ components of $E$ (\cite[Theorem
6.1]{wa:van}).
\end{remark}

In any case we now get from the sequence~\eqref{thetaZ}, the equality
$h^1(\Theta_{\widetilde{X}}(Z))= h^1(\Theta_{\widetilde{X}}) -
\chi(\Theta_{\widetilde{X}}\otimes {\mathcal O}_Z(Z))$. From the
Leray
spectral sequence for $\Theta_{\widetilde{X}}$, we get
$h^1(\Theta_{\widehat{X}}) = h^1(\Theta_{\widetilde{X}}) -
h^0(R^1p_{\star}\Theta_{\widetilde{X}})$. So finally
\begin{equation*}
\begin{split}
h^1(\Theta_{\widehat{X}}(C))& =h^1(\Theta_{\widetilde{X}}(Z)) -
h^0(R^1p_{\star}\Theta_{\widetilde{X}}(Z)) \\
 & = h^1(\Theta_{\widetilde{X}}) - \chi(\Theta_{\widetilde{X}}\otimes
{\mathcal O}_Z(Z)) -  h^0(R^1p_{\star}\Theta_{\widetilde{X}})\\
 & = h^1(\Theta_{\widehat{X}}) - \chi(\Theta_{\widetilde{X}}\otimes
{\mathcal
O}_Z(Z)) \, .
\end{split}
\end{equation*}
Since $p_{\star}{\mathcal O}_{\widetilde X}(kZ)={\mathcal
O}_{\widehat X}(kC)$ and $R^1p_{\star}{\mathcal O}_{\widetilde
X}(kZ)=0$ by the projection formula, we have $p_{\star}{\mathcal
O}_Z(2Z)={\mathcal
O}_C(2C)$ and $R^1p_{\star}{\mathcal O}_Z(2Z)=0$. Thus
$\chi({\mathcal O}_C(2C))=\chi({\mathcal O}_Z(2Z))$.

Putting all this into formula~\eqref{first_h_formula} we get a new
version
\begin{equation}\label{sec_h_formula}
\begin{split}
\dim_{\mathbb C}T^1_X& =
  h^1(\Theta_{\widehat{X}}) - \chi(\Theta_{\widetilde{X}}\otimes
{\mathcal
O}_Z(Z)) + \chi({\mathcal O}_Z(2Z)) + h^0({\mathcal
T}^1_{\widehat{X}}) +
c(X)\\
& = \dim_{\mathbb C}T^1_{\widehat{X}} -
\chi(\Theta_{\widetilde{X}}\otimes
{\mathcal O}_Z(Z)) + \chi({\mathcal O}_Z(2Z))+ c(X)\, .
\end{split}
\end{equation}

\begin{lemma}\label{e-4}
If $X$ is a rational surface singularity, then
$$\chi({\mathcal O}_Z(2Z)) - \chi(\Theta_{\widetilde{X}}\otimes
{\mathcal
O}_Z(Z)) = e-4\, .$$
\end{lemma}
\begin{proof}
We have $\chi({\mathcal O}_Z(2Z))=2Z^2+1=-2e+3$ by Riemann--Roch.

We compute $\chi(\Theta_{\widetilde{X}}\otimes {\mathcal O}_Z(Z))$ in
a
standard manner (see e.g. \cite[page 109]{bk:inf} for another
example). Since $X$ is rational we may
construct a ``computation sequence" $Z_0,\dots ,Z_l=Z$ where
$Z_0=E_{i_0},
Z_1=Z_0+E_{i_1},\dots, Z_l=Z_{l-1}+E_{i_l}$ with the property
$Z_{k-1}\cdot E_{i_k}=1$.
If we tensor the exact sequence
$$0\rightarrow {\mathcal O}_{Z_k}(Z_k) \rightarrow {\mathcal
O}_{Z_{k+1}}(Z_{k+1}) \rightarrow  {\mathcal
O}_{E_{i_{k+1}}}(Z_{k+1}) \rightarrow 0 $$
with $\Theta_{\widetilde{X}}$, we may compute recursively if we know
$\chi(\Theta_{\widetilde{X}} \otimes {\mathcal
O}_{E_{i_{k+1}}}(Z_{k+1}))$. To compute this consider the standard
exact sequence
$$0\rightarrow \Theta_{E_{i_k}} \rightarrow
\Theta_{\widetilde{X}}\otimes
{\mathcal O}_{E_{i_k}} \rightarrow  {\mathcal O}_{E_{i_k}}(E_{i_k})
\rightarrow 0\, .$$ If $k>0$, then after twisting with $Z_k$, we get
this sequence on $E_{i_k}\simeq {\mathbb P}^1$:
$$0\rightarrow {\mathcal O}_{E_{i_k}}(3-b_{i_k}) \rightarrow
\Theta_{\widetilde{X}}\otimes {\mathcal O}_{E_{i_k}}(Z_k)
\rightarrow
{\mathcal O}_{E_{i_k}}(-2b_{i_k}+1) \rightarrow 0$$
where $b_i=-E_i^2$. If $k=0$ subtract $1$ from the degrees of the
left and right sheaves.

After adding everything up we get
$$\chi(\Theta_{\widetilde{X}}\otimes {\mathcal O}_Z(Z)) =
3\sum_{k=0}^l
(2-b_{i_k})-2\, .$$
If $K$ is a canonical divisor, then by the adjunction formula we find
that $-e+1=Z^2=-2-K\cdot Z= -2+\sum_{k=0}^l (2-b_{i_k})$. So
$\chi(\Theta_{\widetilde{X}}\otimes {\mathcal O}_Z(Z)) = -3e+7$ and
we
have proven the lemma.
\end{proof}

\begin{remark} It may be just a curiosity, but the number $e-4$ comes
from sheaves of more deformation theoretical interest. Notice that
$$\chi({\mathcal O}_Z(2Z)) - \chi(\Theta_{\widetilde{X}}\otimes
{\mathcal
O}_Z(Z)) = \chi({\mathcal T}^1_Z(Z)) - \chi(\Theta_Z(Z))\, .$$
This follows from the standard sequence for ${\mathcal T}^1_Z$.
\end{remark}

If we plug the result of Lemma~\ref{e-4} into
formula~\ref{sec_h_formula} we get
\begin{theorem}\label{t1}
If $X$ is a rational surface singularity of embedding dimension $e$
and $\widehat{X}$ is the blow up of $X$, then
$$\dim_{\mathbb C}T^1_{X}=(e-4) + \dim_{\mathbb C}T^1_{\widehat{X}} +
c(X)\, .$$
\end{theorem}

\section{About the correction term $c(X)$.} \label{Cor_term}
\subsection{Alternative definitions}
We
have not been able to compute $c(X)$ in general, though there are
partial results which we present here. First let us list several
other $H^1$´s which have dimension $c(X)$.

\setcounter{temp}{0}
\begin{proposition}\label{alt_def}
If $X$ is a rational surface singularity, then $c(X)$ equals the
dimension of
\begin{list}{\textup{(\roman{temp})}}{\usecounter{temp}}
        \item $H^1(\widehat{X},{\mathfrak m}{\mathcal F}^1/
        {\mathfrak m}^2{\mathcal F}^1)$
        \item $H^1(\widehat{X},N_A({\mathcal O}_{\widehat{X}})(-C))$
        \item $H^1(\widehat{X},{\mathcal F}^1(-C))$
        \item $H^1(\widehat{X},{\mathcal N}_{{\widehat{X}}/M}(C))$.
\end{list}
\end{proposition}
\begin{proof}
To prove (i) it is enough to show that ${\mathfrak
m}^2{\mathcal F}^1$ has support at points. We claim that the
isomorphism in Proposition~\ref{identities} (iii)
induces locally a surjection $T^1_B\twoheadrightarrow x^2T^1_A(B)$.
Indeed if $[\phi]$ is in the kernel $K$ of $T^1_A(B)\to T^1_B$, then
$x^2\phi(f_j)=bx^2\,{\partial g_j}/{\partial x}+\sum
b_ix^2\,{\partial g_j}/{\partial t_i}=b({\partial f_j}/{\partial
x}+\sum t_i\,{\partial f_j}/{\partial x_i})+\sum b_ix\,{\partial
f_j}/{\partial x_i}$ by \eqref{chainrule}. Thus $x^2[\phi]=0$ in
$T^1_A(B)$ and $K$ is contained in the kernel of the multiplication
map $T^1_A(B) \stackrel{\cdot x^2}{\rightarrow} T^1_A(B)$.

Consider the commutative diagram with exact rows and surjective
vertical maps.
\begin{equation*}
\begin{CD}
0 @>>> N_A({\mathcal O}_{\widehat{X}})(-C) @>>> N_A({\mathcal
O}_{\widehat{X}}) @>>> N_A({\mathcal O}_{\widehat{X}})\rvert_C @>>> 0\\
@.  @VV{\alpha}V @VV{\beta}V @VV{\gamma}V \\
0 @>>>{\mathfrak m}{\mathcal F}^1 @>>>{\mathcal F}^1 @>>>{\mathcal
F}^1\rvert_C
@>>> 0
\end{CD}
\end{equation*}
The argument in the proof of Proposition~\ref{h0-zero} shows that
$H^0(N_A({\mathcal O}_{\widehat{X}})\rvert_C)=0$ as well, so
$H^0(\Ker(\gamma))=0$ and $H^1(\Ker(\alpha))$ injects into
$H^1(\Ker(\beta))$. But $\Ker(\beta)$ is an image of
$\pi^{\star}\Der(P)\simeq
e{\mathcal O}_{\widehat{X}}$ so $H^1(\Ker(\beta))=0$. This proves
(ii).

On the other hand $\alpha$ factors surjectively through ${\mathcal
F}^1(-C)$, which gives (iii). We have $H^1({\mathcal
N}_{{\widehat{X}}/M}(C))\simeq
H^1(N_A({\mathcal O}_{\widehat{X}})(-C))$ by
Proposition~\ref{identities}.
\end{proof}

\subsection{Partial results}
\label{section:partial:results}
The exact sequence in Proposition~\ref{identities} sits in the
following large commutative diagram of exact rows and columns.
\begin{equation}
\label{big:diagram}
\begin{CD}
@.  @.  0 @.  0 \\
@.  @.  @VVV @VVV \\
@.  @.  {\mathfrak m}^2{\mathcal F}^1 @>>>
{\mathfrak m}{\mathcal T}^1_{\widehat{X}}(C) @>>> 0\\
@.  @.  @VVV @VVV \\
0 @>>> {\mathcal A}^1_{C/{\widehat{X}}}(C) @>>> {\mathfrak m}{\mathcal F}^1
@>>> {\mathcal T}^1_{\widehat{X}}(C) @>>> 0\\
@.  @VVV @VVV @VVV \\
0 @>>> {\mathcal T}^1_{C\subset\widehat{X}}(C) @>>> {\mathcal T}^1_C(C) @>>>
{\mathcal T}^1_{\widehat{X}}({\mathcal O}_C)(C) @>>> 0\\
@.  @VVV @VVV @VVV \\
@.  0 @>>> {\mathcal T}^2_{\widehat{X}} @>>> {\mathcal T}^2_{\widehat{X}}
@>>> 0
\end{CD}
\end{equation}
We will not prove this, only explain the sheaves and some of the
sequences. We use the notation of Section~\ref{Computation}.

The sheaf ${\mathcal T}^1_{C\subset\widehat{X}}$ is the cokernel of
$\Der({\mathcal O}_{\widehat X},{\mathcal O}_C)\rightarrow {\mathcal
O}_C(C)$ defined locally by $D\mapsto D(x)$. The sheaf ${\mathcal
T}^1_{\widehat{X}}({\mathcal O}_C)$ is locally $T^1_B(B/(x))$. The
right
vertical sequence is induced from the exact sequence
$ 0\rightarrow {\mathcal O}_{\widehat X} \rightarrow {\mathcal
O}_{\widehat
X}(C) \rightarrow {\mathcal O}_C(C) \rightarrow 0$
and the next to bottom row is from the Zariski--Jacobi sequence for
${\mathbb C}\to {\mathcal O}_{\widehat X}\to {\mathcal O}_C$. The map
${\mathfrak m}{\mathcal F}^1 \to {\mathcal T}^1_C(C)$ is locally the
map $xT^1_A(B)\to T^1_{B/(x)}$ which takes $x[\phi]$ to the class of
the homomorphism $\overline{g}_j\mapsto \phi(f_j) + (x)$. Here
$\overline{g}_j$ is the image of $g_j$ in $P_B/(x)$.

With the help of this diagram we can prove
\setcounter{temp}{0}
\begin{proposition}
\label{c_1}
If $X$ is a rational surface singularity, then
\begin{list}{\textup{(\roman{temp})}}{\usecounter{temp}}
        \item $c(X)=0$ if the fundamental cycle is reduced.
        \item $c(X)\ge \dim_{\mathbb C}H^1(\widehat{X},{\mathcal T}^1_C(C))$.
        \item If ${\mathcal T}^2_{\widehat{X}}=0$, i.e. the singularities on
$\widehat{X}$
                all have embedding dimension $\le 4$, then $c(X)=
                \dim_{\mathbb C}H^1(\widehat{X},{\mathcal T}^1_C(C))$.
        \end{list}
\end{proposition}
\begin{proof}
If the fundamental divisor is reduced, then $C$ is reduced
(\cite[Proof of Corollary 3.6]{wa:equ}). In this case ${\mathcal
A}^1_{C/{\widehat{X}}}$ has support at points so $H^1({\mathcal
A}^1_{C/{\widehat{X}}}(C))=0$.

From the diagram we get an exact sequence
$ 0\rightarrow {\mathfrak m}{\mathcal F}^1/{\mathfrak m}^2{\mathcal
F}^1 \rightarrow {\mathcal T}^1_C(C) \rightarrow {\mathcal
T}^2_{\widehat{X}}$, and $c(X)=h^1({\mathfrak m}{\mathcal
F}^1/{\mathfrak m}^2{\mathcal F}^1)$ by Proposition~\ref{alt_def}.
This proves the last two statements.
\end{proof}

Proposition~\ref{c_1} allows us to generate examples where $c(X)>0$.
\begin{proposition}
\label{c_2}
If $X$ is a rational surface singularity, $\widehat{X}$ is smooth and
$C$
is non-reduced, then $c(X)>0$.
\end{proposition}
\begin{proof}
In this case $\widehat{X}=\widetilde{X}$, $C=Z$ is the fundamental
divisor
and $c(X)=h^1({\mathcal T}^1_Z(Z))$. It follows from
\cite[2.6]{wa:sim}
that $H^1({\mathcal T}^1_Z(Z))\simeq H^1({\mathcal O}_{Z-E}(2Z))$
where $E=Z_{\text{red}}$.

We claim that $H^0({\mathcal O}_{Z-E}(2Z))=0$. In \cite{wa:van} Wahl
proves that $H^1_E({\mathcal O}_{\widetilde{X}}(E))=0$. The proof
actually shows that $H^1_E({\mathcal
O}_{\widetilde{X}}(E+Z^\prime))=0$
for any cycle $Z^\prime$ with $Z^\prime\cdot E_i\le 0$ for all
irreducible components $E_i$ of $E$. In particular $H^1_E({\mathcal
O}_{\widetilde{X}}(E+Z))=0$, and we may use \cite[Lemma B.2]{wa:equi}
to
conclude that $H^0({\mathcal O}_Z(E+2Z))=0$. But
$H^0({\mathcal O}_{Z-E}(2Z))$ injects into this last $H^0$.

Thus $h^1({\mathcal T}^1_Z(Z))=-\chi({\mathcal
O}_{Z-E}(2Z))=(Z-E)\cdot (K-Z)$ using Riemann--Roch. (See \cite[Proof
of Proposition 2.15]{wa:sim} for a similar argument.) If we set
$b_i=-E_i^2$, $Z=\sum n_iE_i$ and $r_i=-Z\cdot E_i$ we find
$(Z-E)\cdot (K-Z)=\sum (n_i-1)(b_i-2+r_i)$. Now we have assumed
$\widehat{X}=\widetilde{X}$, so all the $r_i>0$ by Theorem~\ref{Tyu}.
Thus
$c(X)>0$ in this case if $Z\ne E$; i.e. at least one $n_i \ge 2$.
\end{proof}

It is a purely combinatorial problem to make dual graphs for
rational singularities satisfying the conditions in
Proposition~\ref{c_2}. The one with lowest multiplicity is the
``standard counter example" to the $T^1$ and $T^2$ formulas appearing
before this paper -- see e.g. \cite{bk:inf}. Here is the dual graph:
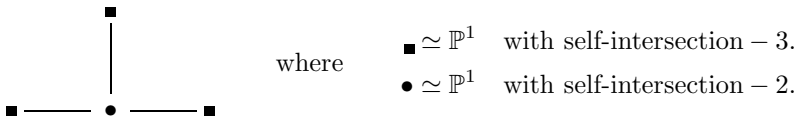
\begin{equation*}
\begin{aligned}
\begin{diagram}
        \node[2]{\makebox(0,0){\rule{1.2mm}{1.2mm}}} \arrow{s,-}\\
    \node{\makebox(0,0){\rule{1.2mm}{1.2mm}}} \arrow{e,-}
\node{\bullet}
                \arrow{e,-} \node{\makebox(0,0){\rule{1.2mm}{1.2mm}}}
\end{diagram}
\end{aligned}
\qquad\text{where}\qquad
\begin{aligned}
{\makebox(0,0){\rule{1.2mm}{1.2mm}}}\, &\simeq {\mathbb P}^1
        \quad\text{with self-intersection} -3.\\
\bullet &\simeq {\mathbb P}^1
        \quad\text{with self-intersection} -2.
\end{aligned}
\end{equation*}
In fact any exceptional configuration of $4$ components with this
type of intersection will have $c(X)>0$, as long as the central curve
has self-intersection $-2$ and the other self-intersections are $\le
-3$.

If one extends the three arms off the central $-2$ curve, then these
singularities will also have $c(X)>0$ as long as neighbors of the
$-2$ curve have self-intersection $\le -4$ if the arm has length $>
1$ and non-end nodes have self-intersection $\le -3$.

Here is an example with $e=7$.
\begin{equation*}
\begin{aligned}
\begin{diagram}
        \node[2]{\makebox(0,0){\rule{1.2mm}{1.2mm}}} \arrow{s,-}\\
    \node{\makebox(0,0){\rule{1.2mm}{1.2mm}}} \arrow{e,-}
\node{\bullet}
                \arrow{e,-} \node{\times} \arrow{e,-}
                \node{\bullet}
\end{diagram}
\end{aligned}
\qquad\text{where}\qquad
\begin{aligned}
\times &\simeq {\mathbb P}^1
        \quad\text{with self-intersection} -4.\\
{\makebox(0,0){\rule{1.2mm}{1.2mm}}}\, &\simeq {\mathbb P}^1
        \quad\text{with self-intersection} -3.\\
\bullet &\simeq {\mathbb P}^1
        \quad\text{with self-intersection} -2.
\end{aligned}
\end{equation*}

\bibliographystyle{amsalpha}
\bibliography{dr}
\end{document}